\documentclass[a4paper,12pt]{article}
\usepackage{amsmath,amsthm,amsfonts,amssymb,graphicx,url,color}

\date{}

\newcommand{\sinc}{\operatorname{sinc}}
\theoremstyle{theorem}
\newtheorem{theorem}{Theorem}

\theoremstyle{definition}

\newtheorem*{remark}{Remark}

\newtheorem{lemma}[theorem]{Lemma}
\newtheorem{proposition}[theorem]{Proposition}

\newenvironment{example}[1][Example]{\begin{trivlist}
\item[\hskip \labelsep {\bfseries #1}]}{\end{trivlist}}

\begin{document}

\title{ More remarkable sinc integrals and sums }

\author{Gert Almkvist and Jan Gustavsson}

\maketitle

\bigskip

\begin{abstract}
We use Poisson summation formula to calculate integrals of products of sinc functions (cf. \cite{4}) and related integrals as in \cite{5} and \cite{3}. 
We also generalize the one in \cite{5} and introduce other remarkable integrals.

Finally we give a sum version of Siegel-type lower bound. (cf. \cite{2},  Theorem 3)


\end{abstract}

\bigskip

\subsection*{Introduction.}

\bigskip

In 2001 David and Jonathan Borwein in \cite{4} via Fourier transform theory proved that
\begin{equation}
\int_{0}^{\infty}\prod_{k=0}^{n}\sinc(\dfrac{t}{2k+1})\, dt = \dfrac{\pi}{2}, \quad n = 0, 1, 2, 3, 4, 5, 6
\end{equation}
but less than $\dfrac{\pi}{2}$ if $n>6$. Quite recently the integral came up again in \cite{5} where the Fourier transformation is clarified by a graphic approach. H. Schmid also proves that
\begin{equation}
\int_{0}^{\infty}2\cos(t)\prod_{k=0}^{n}\sinc(\dfrac{t}{2k+1})\, dt = \dfrac{\pi}{2}, \quad n = 0, 1, 2, \dots , 55
\end{equation}
but less than $\dfrac{\pi}{2}$ if $n>55$.
We will use the Poisson summation formula to prove both (1) and (2). The formula gives us also a generalization of (2) and other curious sinc integrals. E.g. we find that
\begin{equation}
\int_{0}^{\infty}(2\cos(t)+2\cos(3t))\prod_{k=0}^{n}\sinc(\dfrac{t}{2k+1})\, dt = \dfrac{\pi}{2}, \quad n = 0, 1, 2, \dots , 3090.
\end{equation}
For $n> 3090$ the value of the integral will be $< \dfrac{\pi}{2}$.

The other day B\"{a}sel in \cite{3} remarked that (2) via an elementary formula can be deduced from (1).

In \cite{2} Baillie, D. Borwein and J. M. Borwein prove that (Theorem 3)
\begin{equation*}
\int_{0}^{\infty}\prod_{k=0}^{n}\sinc(a_k x)\, dx \ge \int_{0}^{\infty}\sinc^{n+1}(a_0 x)\, dx
\end{equation*}
where
\begin{equation*}
a_0 \ge a_k >0 \text{ for } k= 1, 2, \ldots , n.
\end{equation*}
They also write: ''Perhaps a somewhat analogous version of Theorem 3 holds for sums?'' We will give an answer.

\bigskip

\subsection*{Theoretical tools.}

\bigskip

For appropriate (see \cite{6}) functions the Poisson summation formula may be stated as
\begin{equation}
\sum_{k=-\infty}^{\infty}\hat{f}(\omega - k\Omega) = T\sum_{m=-\infty}^{\infty}f(mT)e^{-im\omega T}
\end{equation}
where $T\Omega = 2\pi$ and $\displaystyle \hat{f}(\omega) = \int_{-\infty}^{\infty}e^{-i\omega t}f(t)\, dt$ (Fourier transform).

If we choose $T = 1$ and $\omega = 0$ and if we also assume $f$ to be an even function then we can write (4) as
\begin{equation}
\sum_{m=-\infty}^{\infty}f(m) = \hat{f}(0) + 2\sum_{k=1}^{\infty}\hat{f}(2k\pi) = \int_{-\infty}^{\infty}f(t)\, dt + 2\sum_{k=1}^{\infty}\hat{f}(2k\pi).
\end{equation}
If we instead choose $\omega  = \pi$ we can write (4) as
\begin{gather}
\sum_{m=-\infty}^{\infty}f(m)e^{-im\pi} = \hat{f}(\pi)+\hat{f}(-\pi) + 2\sum_{k=1}^{\infty}\hat{f}(\pi +2k\pi) \notag\\ = \int_{-\infty}^{\infty}2\cos(\pi t)f(t)\, dt + 2\sum_{k=1}^{\infty}\hat{f}(\pi +2k\pi).
\end{gather}
We sum this up as
\begin{lemma}
\textit{Let} $f(t)$ \textit{be even and sufficiently summable and integrable. Then}
\begin{equation}
\sum_{m=-\infty}^{\infty}f(m) = \int_{-\infty}^{\infty}f(t)\, dt \Leftrightarrow \sum_{k=1}^{\infty}\hat{f}(2k\pi) = 0
\end{equation}
and
\begin{equation}
\sum_{m=-\infty}^{\infty}f(m)e^{-im\pi} = \int_{-\infty}^{\infty}2\cos(\pi t)f(t)\, dt \Leftrightarrow 
\sum_{k=1}^{\infty}\hat{f}(\pi +2k\pi) = 0.
\end{equation}
\end{lemma}
From this we build 
\setcounter{theorem}{0}
\begin{theorem}
Assume that $f_0, f_1, \ldots , f_n$ are even functions and assume that their Fourier transforms have their supports in $[-1,1]$. Then
\begin{equation}
\sum_{m=-\infty}^{\infty}f_0(a_0m)f_1(a_1m)\cdot\ldots\cdot f_n(a_nm) = \int_{-\infty}^{\infty}f_0(a_0t)f_1(a_1t)\cdot\ldots\cdot f_n(a_nt)\, dt
\end{equation}
if $a_0, a_1,\ldots , a_n$ are positive and
\begin{equation}
\sum_{k=0}^{n}a_k < 2\pi.
\end{equation}
Furthermore
\begin{gather}
\sum_{m=-\infty}^{\infty}f_0(a_0m)f_1(a_1m)\cdot\ldots\cdot  f_n(a_nm)e^{-im\pi} \notag\\= \int_{-\infty}^{\infty}2\cos(\pi t)f_0(a_0t)f_1(a_1t)\cdot\ldots\cdot f_n(a_nt)\, dt
\end{gather}
if $a_0, a_1,\ldots , a_n$ are positive and
\begin{equation}
\sum_{k=0}^{n}a_k < 3\pi.
\end{equation}
\end{theorem}

\bigskip

\begin{proof}
If $f(t)$ is a function with a Fourier transform with support in $[-1,1]$ then $f(at)$ has a Fourier transform with support in $[-a,a]$. Put $f(t) = f_1(a_1t)f_2(a_2t)\cdot\ldots\cdot f_n(a_nt)$. The Fourier transform $\hat{f}$ is given by a convolution with support in $[-(a_0+a_1+\ldots +a_n), a_0+a_1+\ldots +a_n]$. 
Thus (9) is a consequence of (7).

In a similar way we get (11) from (8).
\end{proof}
\bigskip

\begin{remark}
We used (7) and (9) already in \cite{1}.
\end{remark}

\bigskip

The idea of the next proposition can be found in \cite{4}.

\setcounter{theorem}{0}

\begin{proposition}
If $f(t)$ is a function in $L_2(\mathbf{R})$ with a Fourier transform with support in $[-a,a]$ and if $f$ is continuous at $t = 0$ then
\begin{equation}
\int_{-\infty}^{\infty}f(t)\frac{\sin(bt)}{t}\, dt = \pi f(0)
\end{equation}
if $0<a<b$.
\end{proposition}

\bigskip

\begin{proof}
We observe that $\dfrac{\sin(bt)}{t}$ has the Fourier transform $\pi (H(\omega +b) - H(\omega -b))$, where $H$ is the Heaviside function given by
\begin{equation*}
H(\omega) = \left\{\begin{array}{ll}
1, & \text{if } \omega > 0\\
0,& \text{if } \omega < 0.
\end{array}\right.
\end{equation*}
Put $f_{\varepsilon}(t) = f(t)\sinc(\varepsilon t)$ where $a<a+\varepsilon < b$. Then $f_{\varepsilon}$ is in $L_1(\mathbf{R})\cap L_2(\mathbf{R})$. According to Parseval's theorem $\hat{f_{\varepsilon}}$ has its support in $[-a-\varepsilon  , a+\varepsilon]$. 

A version of Fourier inversion formula states that
\begin{equation*}
g(t) = \dfrac{1}{2\pi}\int_{\infty}^{\infty}e^{i\omega t}\hat{g}(\omega)\, d\omega
\end{equation*}
if $g$ and $\hat{g}$ are in $L_1(\mathbf{R})$ and $g$ is continuous at $t$.

If we combine Parseval's theorem and the Fourier inversion formula we get
\begin{gather*}
\int_{-\infty}^{\infty}f_{\varepsilon}(t)\frac{\sin(bt)}{t}\, dt = \dfrac{1}{2\pi}\int_{-\infty}^{\infty}\hat{f_{\varepsilon}}(\omega)\pi (H(\omega +b) - H(\omega -b))\, d\omega \\= \dfrac{\pi}{2\pi}\int_{-a-\varepsilon}^{a+\varepsilon}\hat{f_{\varepsilon}}(\omega)\, d\omega = \pi\dfrac{1}{2\pi}\int_{-\infty}^{\infty}e^{i0\omega}\hat{f_{\varepsilon}}(\omega)\, d\omega = \pi f_{\varepsilon}(0) = \pi f(0).
\end{gather*}
But since
\begin{equation*}
\left|f_{\varepsilon}(t)\frac{\sin(bt)}{t}\right| \le \left|f(t)\frac{\sin(bt)}{t}\right| \in L_1(\mathbf{R})
\end{equation*}
we can use Lebesgue's dominated convergence theorem and let $\varepsilon \to  0$. This will give us (13).
\end{proof}

\bigskip

\subsection*{Applications or examples.}

\bigskip

\begin{example}\textbf{1.}
We will here study (1). If we change the variable $t$ to $\pi t$ and observe that the integrand is even we have to prove that
\begin{equation}
\int_{-\infty}^{\infty}\prod_{k=0}^{n}\sinc(\dfrac{\pi t}{2k+1})\, dt = 1, \quad n = 0, 1, 2, 3, 4, 5, 6
\end{equation}
and $< 1$ if $n>6.$

The function $\sinc(t)$ has the Fourier transform $\pi (H(\omega +1)-H(\omega -1))$. Let the functions $f_0, f_1, \ldots , f_n$ in (9) all be $\sinc(t)$. Consequently $a_k = \dfrac{\pi}{2k+1}$. The condition (10) corresponds to $\displaystyle \sum_{k=0}^{n}\dfrac{1}{2k+1} < 2.$ Since
\begin{equation*}
\sum_{k=0}^{6}\dfrac{1}{2k+1} = \dfrac{88069}{45045} \approx 1.9551
\end{equation*}
but
\begin{equation*}
\sum_{k=0}^{7}\dfrac{1}{2k+1} = \dfrac{91072}{45045} \approx 2.0218
\end{equation*}
(10) is fulfilled for $n= 0, 1, \ldots , 6$ but not for $n >6.$

The value of the integral in (14) will now be given by the series on the left hand side of (9). Since
\begin{equation*}
f_0(a_0m)= \sinc(\pi m) = \left\{\begin{array}{ll}
1,&\text{if } m = 0\\
0,&\text{otherwise}
\end{array}\right.
\end{equation*}
the sum of the series will be $1$.

If $n>6$ the integral in (14) will be
\begin{equation*}
1- 2\sum_{k=1}^{\infty}\hat{f}(2k\pi) < 1. 
\end{equation*}
Cf. (5) and remember that $\hat{f}(\omega)$ is a convolution which is positive for $\omega \in [-(a_0+a_1+ \ldots + a_n), a_0+a_1+ \ldots + a_n]$
\end{example}

\bigskip

\begin{example}\textbf{2.}
Here we prove (2) and calculate the integral when $n=56.$ It is difficult to test the calculation with a computer. A straightforward calculation of the integral by means of Maple was no success.

We copy the method in Example 1. Via (11) we find that
\begin{equation*}
\int_{-\infty}^{\infty}2\cos(\pi t)\prod_{k=0}^{n}\sinc(\dfrac{\pi t}{2k+1})\, dt = 1
\end{equation*}
if \, $\displaystyle \sum_{k=0}^{n}\dfrac{1}{2k+1} < 3.$ But
\begin{equation*}
\sum_{k=0}^{55}\dfrac{1}{2k+1} \approx 2.994437501
\end{equation*}
and
\begin{equation*}
\sum_{k=0}^{56}\dfrac{1}{2k+1} \approx 3.003287059.
\end{equation*}
and we have proved (2).

In order to calculate
\begin{equation*}
\int_{-\infty}^{\infty}2\cos(\pi t)\prod_{k=0}^{56}\sinc(\dfrac{\pi t}{2k+1})\, dt
\end{equation*}
we use (6) with 
\begin{equation*}
f(t) = \prod_{k=0}^{56}\sinc(\dfrac{\pi t}{2k+1}).
\end{equation*}
Since $f(m) = 0$ for $m \neq 0$ and $f(0) = 1$ we get that
\begin{equation*}
\int_{-\infty}^{\infty}2\cos(\pi t)\prod_{k=0}^{56}\sinc(\dfrac{\pi t}{2k+1})\, dt = 1 - 2\sum_{k=1}^{\infty}\hat{f}({\pi +2k\pi}).
\end{equation*}
But the support of $\hat{f}$ is $[-a,a]$, where
\begin{equation*}
a = \sum_{k=0}^{56}\dfrac{\pi}{2k+1}\approx 9.435104562.
\end{equation*}
Thus $3\pi < a < 5\pi$ and 
\begin{equation*}
\int_{-\infty}^{\infty}2\cos(\pi t)\prod_{k=0}^{56}\sinc(\dfrac{\pi t}{2k+1})\, dt = 1 - 2\hat{f}(3\pi).
\end{equation*}
We have to investigate the Fourier transform $\hat{f}$, which we also denote by $F$.  We can rewrite $f(t)$ as
\begin{equation*}
f(t) = \prod_{k=0}^{56}\sinc(\dfrac{\pi t}{2k+1}) = \dfrac{c}{t^{57}{\pi}^{57}}\prod_{k=0}^{56}\sin(\dfrac{\pi t}{2k+1}) 
\end{equation*}
where $\displaystyle c = \prod_{k=1}^{56}(2k+1)$. Via Euler's formulae we continue to
\begin{equation}
t^{57}f(t) = \dfrac{c}{{(2\pi)}^{57}i}\prod_{k=0}^{56}\left(e^{i\frac{\pi t}{2k+1}}-e^{-i\frac{\pi t}{2k+1}}\right) = \dfrac{c}{{(2\pi)}^{57}i}\sum_{k}e_ke^{id_kt}
\end{equation}
where $e_k$ are coefficients and the $d_k$ are numbers between $-a$ and $a$. The largest $d_k$ is $a$ and the next largest is $\displaystyle \sum_{k=0}^{55}\frac{\pi}{2k+1} - \dfrac{\pi}{113} \approx 9.379501153$. In the interval between that number and $a$ we find $3\pi$.  On that interval $F(\omega)$ will be a polynomial.

The Fourier transform of $tf(t)$ is $iF'(\omega)$ (at least in the distributional sense) and the Fourier transform of $t^{57}f(t)$ is $i^{57}F^{(57)}(\omega)
 = iF^{(57)}(\omega)$. Thus if we apply the Fourier transformation to both sides in (15) we get
\begin{equation*}
iF^{(57)}({\omega}) = \dfrac{c}{{(2\pi)}^{57}i}(2\pi\delta(\omega -a) + \text{other } \delta \text{terms })
\end{equation*}
where $\delta$ is the Dirac ''function''. Notice that the Fourier transform of $\displaystyle e^{id_kt}$ is $\displaystyle 2\pi\delta(\omega -d_k)$.

From \begin{equation*}
F^{(57)}(\omega) = - \dfrac{c}{{(2\pi)}^{56}}(\delta(\omega -a) + \text{other } \delta \text{terms })
\end{equation*}
we conclude that
\begin{equation*}
F^{(56)}(\omega) =  \dfrac{c}{{(2\pi)}^{56}} \quad \text{ if } \quad \sum_{k=0}^{55}\frac{\pi}{2k+1} - \dfrac{\pi}{113} < \omega < a.
\end{equation*}
$F$ is built up of 57 convolution factors. Thus $F$ is 55 times differentiable. But $F(\omega) \equiv 0$ if $\omega > a$. All together gives us that
\begin{equation}
 F(\omega) = \dfrac{c}{2^{56}{\pi}^{56}56!}(\omega -a)^{56} \quad \text{ if } \quad \sum_{k=0}^{55}\frac{\pi}{2k+1} - \dfrac{\pi}{113} < \omega < a.
 \end{equation}
Introduce $\displaystyle b = \sum_{k=0}^{56}\dfrac{1}{2k+1}$. Thus $a = b\pi$ and we can write (16) as
\begin{equation*}
F(\omega) = \dfrac{c}{2^{56}{\pi}^{56}56!}(\omega -b\pi)^{56}.
\end{equation*}
Now we get
\begin{equation*}
2F(3\pi) = \dfrac{c}{2^{55}{\pi}^{56}56!}(3\pi -b\pi)^{56} = \dfrac{c}{2^{55}56!}(3 -b)^{56} = \dfrac{N}{D}
\end{equation*}
where
\begin{equation*}
N = 347^{56}\cdot 39608671351^{56}\cdot 1786013712647720237751897933348037^{56}
\end{equation*}
and
\begin{gather*}
D = 2^{53}\cdot 3^{222}\cdot 5^{112}\cdot 13^{56}\cdot 41^{56}\cdot 59^{55}\cdot 61^{55}\cdot 67^{55}\cdot 71^{55}\cdot 73^{55}\cdot 79^{55}\cdot 11^{56}\cdot 17^{56}\\\cdot 7^{112}\cdot 47^{56}\cdot 29^{55}\cdot 31^{55}\cdot 37^{55}\cdot 43^{56}\cdot 53^{56}\cdot 23^{56}\cdot 19^{55}\cdot 83^{55}\cdot 89^{55}\cdot 97^{55}\cdot 101^{55}\\\cdot 103^{55}\cdot 107^{55}\cdot 109^{55}\cdot 113^{55}.
\end{gather*}
Finally we get
\begin{gather*}
\int_{-\infty}^{\infty}2\cos(\pi t)\prod_{k=0}^{56}\sinc(\dfrac{\pi t}{2k+1})\, dt = 1 - 2F(3\pi) = 1 - \dfrac{N}{D}\\ \approx 1 - 1.484870809\cdot 10^{-138}.
\end{gather*}
The integral is apparently only slightly smaller than 1.
\end{example}

\bigskip

\begin{example}\textbf{3.}
In order to prove (3) we prove that
\begin{equation*}
\int_{-\infty}^{\infty}(2\cos(\pi t)+2\cos(3\pi t))\prod_{k=0}^{n}\sinc(\dfrac{\pi t}{2k+1})\, dt = 1, \quad n = 0, 1, 2, \dots , 3090.
\end{equation*}
We write (6) as
\begin{equation*}
\sum_{m=-\infty}^{\infty}f(m)e^{-im\pi} = \int_{-\infty}^{\infty}(2\cos(\pi t) + 2\cos(3\pi t))f(t)\, dt + 2\sum_{k=2}^{\infty}\hat{f}(\pi +2k\pi).
\end{equation*}
To finish the argument we observe that
\begin{equation*}
\sum_{k=0}^{3090}\dfrac{\pi}{2k+1} \approx 15.70758624 < 5\pi \approx 15.70796327 < \sum_{k=0}^{3091}\dfrac{\pi}{2k+1} \approx 15.70809434.
\end{equation*}
With the same technique we get
\begin{gather*}
\int_{-\infty}^{\infty}(2\cos(\pi t)+2\cos(3\pi t) + 2\cos(5\pi t))\prod_{k=0}^{n}\sinc(\dfrac{\pi t}{2k+1})\, dt = 1, \\ n = 0, 1, 2, \dots , 168802.
\end{gather*}
\end{example}

\bigskip

\begin{example}\textbf{4.}
It is not quite easy to find a pattern for  ''the breaking points'' $n = 6, 3090$ and $168802$, which we have met above. But if we change $\displaystyle \prod_{k=0}^{n}\sinc(\dfrac{\pi t}{2k+1})$ to $\displaystyle \sinc^{n}(\pi t)$ this problem can be solved.

With experience from the other examples we now rewrite (6) as
\begin{equation}
\sum_{m=-\infty}^{\infty}f(m)e^{-im\pi} =  \int_{-\infty}^{\infty}2\sum_{k=0}^{m}\cos((2k+1)\pi t)f(t)\, dt + 2\sum_{k=m+1}^{\infty}\hat{f}((2k+1)\pi).
\end{equation}
Put $f(t) = \displaystyle \sinc^{n}(\pi t)$. Then the left hand side in (17) will be 1 and the points $(2k+1)\pi ,\, k= m+1, m+2, \dots$ will all be outside the support of $\hat{f}$ if $n\pi \le (2m+3)\pi$. Thus
\begin{equation*}
\int_{-\infty}^{\infty}2\sum_{k=0}^{m}\cos((2k+1)\pi t)\sinc^{n}(\pi t)\, dt = \left\{\begin{array}{ll}
1,& \text{if } n = 0, 1, 2, \ldots ,2m+3\\
< 1,& \text{if } n = 2m+4, 2m+5, \ldots 
\end{array}\right.
\end{equation*}
\end{example}
\bigskip

\begin{example}\textbf{5.}
Let us introduce a function which is not a sinc function. Put \begin{equation*}
f(t) = \dfrac{t\sin(t)-\cos(t) + e}{(1+t^{2})(e-1)}.
\end{equation*}
Then $f(0) = 1$ and
\begin{equation*}
\hat{f}(\omega) = \left\{\begin{array}{ll}
\dfrac{\pi}{1-e^{-1}}e^{-|\omega|},&\text{if } |\omega| < 1\\[2ex]
0, \text{otherwise.}
\end{array}\right.
\end{equation*}
Now we get another curious integral
\begin{equation*}
\int_{-\infty}^{\infty}f(a_0t)f(a_1t)\cdot\ldots\cdot f(a_nt)\dfrac{\sin(bt)}{t}\, dt = \pi, \quad \text{if } \, a_0 + a_1 + \dots + a_n < b.
\end{equation*}
We also assume that all $a_k >0$ and that $b>0.$
\begin{proof}
\textbf{Alternative 1.} Proposition 1 delivers a proof directly.

\bigskip

\textbf{Alternative 2.} If we change $t$ to $\dfrac{\pi}{b}t$ we get
\begin{gather*}
\int_{-\infty}^{\infty}f(a_0t)f(a_1t)\cdot\ldots\cdot   f(a_nt)\dfrac{\sin(bt)}{t}\, dt\\ = \int_{-\infty}^{\infty}f(\dfrac{\pi a_0}{b}t)f(\dfrac{\pi a_1}{b}t)\cdot\ldots\cdot   f(\dfrac{\pi a_n}{b}t)\dfrac{\sin(\pi t)}{t}\, dt = 1
\end{gather*}
if (see (10))
\begin{equation*}
\dfrac{\pi a_0}{b} + \dfrac{\pi a_1}{b}+ \ldots + \dfrac{\pi a_n}{b} + \pi < 2\pi \Leftrightarrow a_0 + a_ 1 + \ldots + a_n < b.
\end{equation*}
\end{proof}
\end{example}

\bigskip

\begin{remark}
With the same technique as in alternative 2 we can use Theorem 1 to handle the last curious integrals $I_n(b)$ in \cite{3}.
\end{remark}

\bigskip

\begin{example}\textbf{6.}
In \cite{2} a ''Lower Bound'' result is proved. If $a_0 \ge a_k$ for $k = 1, 2, \ldots , n$ then
\begin{equation}
\int_{0}^{\infty}\prod_{k=0}^{n}\sinc(a_k x)\, dx \ge \int_{0}^{\infty}\sinc^{n+1}(a_0 x)\, dx .
\end{equation}
If $(n+1)a_0 < 2\pi$ we can use Theorem 1 to translate the two integrals to sums. Thus
\begin{equation*}
\sum_{m=0}^{\infty}\prod_{k=0}^{n}\sinc(a_k m) \ge \sum_{m=0}^{\infty}\sinc^{n+1}(a_0m).
\end{equation*}
which is an analogous version to (18) for sums.

Since
\begin{equation*}
\sum_{m=0}^{\infty}\sinc(\dfrac{5\pi}{4} m)\sinc^{2}(m) \approx .8999999997 < \sum_{m=0}^{\infty}\sinc^{3}(\dfrac{5\pi}{4} m) \approx .9960000000
\end{equation*}
the condition $(n+1)a_0 < 2\pi$ cannot be omitted.
\end{example}

\vspace*{1cm}

Center for Mathematical Sciences, Lund University, Box 118, SE-22100 Lund, Sweden

\bigskip

gert.almkvist@yahoo.se

jan.gustavsson@math.lth.se or gustavsson.jan@telia.com

\end{document}